\documentclass{amsart}

 \usepackage{amssymb}

\usepackage{amsthm}

\newtheorem{theorem}{Theorem}[section]

\newtheorem{remark}{Remark}[section]

\numberwithin{equation}{section}

\newtheorem*{theorem*}{Theorem}

\newcommand{\R}{\mathbb{R}}
\newcommand{\N}{\mathbb{N}}
\newcommand{\be}{\begin{equation}}
\newcommand{\ee}{\end{equation}}

\begin{document}


\title[Improved Caffarelli-Kohn-Nirenberg and trace inequalities\dots]{Improved Caffarelli-Kohn-Nirenberg and trace inequalities for radial functions}

\author{Pablo L.  De N\'apoli}
\address{Departamento de Matem\'atica \\
Facultad de Ciencias Exactas y Naturales \\ 
Universidad de Buenos Aires \\
Ciudad Universitaria \\
1428 Buenos Aires, Argentina}
\email{pdenapo@dm.uba.ar}

\author{Irene Drelichman}
\address{Departamento de Matem\'atica \\
Facultad de Ciencias Exactas y Naturales \\ 
Universidad de Buenos Aires \\
Ciudad Universitaria \\
1428 Buenos Aires, Argentina}
\email{irene@drelichman.com}

\author{Ricardo G. Dur\'an}
\address{Departamento de Matem\'atica \\
Facultad de Ciencias Exactas y Naturales \\ 
Universidad de Buenos Aires \\
Ciudad Universitaria \\
1428 Buenos Aires, Argentina}
\email{rduran@dm.uba.ar}

\thanks{Supported by ANPCyT under grant PICT 01307 and by Universidad de
Buenos Aires under grants X070 and X837. The first and third authors are members of
CONICET,
Argentina.}

\begin{abstract}
We show that Caffarelli-Kohn-Nirenberg first order interpolation inequalities as well as weighted trace inequalities in $\mathbb{R}^n \times \mathbb{R}_+$ admit a better  range of power weights if we restrict ourselves to the space of radially symmetric functions. 
\end{abstract}


\maketitle
\section{Introduction}

The aim of this paper is to show that inequalities of Caffarelli-Kohn-Nirenberg type hold for a wider class of exponents if we restrict ourselves to the space of radially symmetric functions. To make this precise, recall the classical first-order interpolation inequality obtained in  \cite{CKN}:

\begin{theorem*}[Caffarelli-Kohn-Nirenberg]
Assume 
\begin{equation}
\label{cond1}
p, q \ge 1, \quad r>0, \quad 0\le a \le 1
\end{equation}
\begin{equation}
\label{cond2}
\frac{1}{p}+\frac{\alpha}{n}, \quad \frac{1}{q}+\frac{\beta}{n}, \quad \frac{1}{r}+\frac{\gamma}{n} >0,
\end{equation}
where
\begin{equation}
\label{sigma}
\gamma = a\sigma + (1-a)\beta.
\end{equation}
Then, there exists a positive constant $C$ such that the following inequality holds for all $u\in C_0^\infty(\mathbb{R}^n)$
\begin{equation}
\label{ckn}
\| |x|^\gamma u\|_{L^r} \le C \| |x|^\alpha \nabla u \|_{L^p}^a \| |x|^\beta u\|_{L^q}^{1-a}
\end{equation}
if and only if the following relations hold:
\begin{equation}
\label{reescale}
\frac{1}{r}+\frac{\gamma}{n} =  a \left( \frac{1}{p}+ \frac{\alpha-1}{n}\right) + (1-a)\left( \frac{1}{q}+ \frac{\beta}{n}\right)
\end{equation}
\begin{align}
&0 \le \alpha-\sigma& \qquad &\mbox{ if } \, a>0, &
\\ \nonumber \mbox{and} 
\\ &\alpha-\sigma \le 1&  \qquad &\mbox{ if } \, a>0 \quad \mbox{ and } \quad \frac{1}{p}+\frac{\alpha-1}{n}=\frac{1}{r}+\frac{\gamma}{n}.
\end{align}
\end{theorem*}

Although the conditions of the above theorem cannot be improved in general, we will prove that if  we require $u$ to be radially symmetric, inequality \eqref{ckn} holds true for certain negative values of $\alpha-\sigma$ also. Indeed, the following improvement  holds in this particular case:

\begin{theorem}
\label{ckn-radial}
Assume conditions \eqref{cond1}, \eqref{cond2}, \eqref{sigma} and \eqref{reescale} hold. Then there exists a positive constant $C$ such that inequality \eqref{ckn} holds for all radially symmetric $u\in C_0^\infty(\mathbb{R}^n)$ and all
\begin{equation}
\label{cond-qr}
\frac{1-a}{q} \le \frac{1}{r} \le \frac{a}{p} + \frac{1-a}{q}
\end{equation}
provided that, if $a>0$,
\begin{equation}
\label{alfa-sigma}
(n-1)\left[ \frac{1}{a} \left( \frac{1}{r}-\frac{1}{q}\right)+ \frac{1}{q} - \frac{1}{p}\right] \le \alpha - \sigma \le 0
\end{equation}
and
\begin{equation}
\label{cond-sigma}
-\frac{\sigma}{n} < \frac{1}{a} \left(\frac{1}{r}-\frac{1}{q} \right) + \frac{1}{q},
\end{equation}
with strict inequality in \eqref{alfa-sigma} if $p=1$.
\end{theorem}

\begin{remark}
If $\sigma>0$ condition \eqref{cond-sigma} trivially holds because of \eqref{cond-qr}, and thus our result admits a simpler statement in this case. 
\end{remark}

The key to our proof is to use the well-known inequality relating $u\in C_0^\infty(\mathbb{R}^n)$ with the fractional integral (also called Riesz potential) of its gradient, namely
\begin{equation}
|u(x)| \le C \int_{\mathbb{R}^n} \frac{|\nabla u|(y)}{|x-y|^{n-1}} \, dy =: T_{n-1} (|\nabla u|)(x)
\label{representation}
\end{equation}
together with  improved weighted estimates for fractional integrals of radial functions from \cite{DDD} and the observation that   inequality \eqref{ckn} enjoys a certain self-improving property.  It is worth noting that this method of proof is different from that of the original proof by L. Caffarelli, R. Kohn and L. Nirenberg \cite{CKN}, and also from a different approach developed by  F. Catrina and Z-Q. Wang in \cite{CW}.

We then show that also improved trace inequalities can be obtained in a similar way, but with a slightly different operator involved in the formula \eqref{representation}, for which we prove the required weighted estimates, that play the same role as that played by the result  of \cite{DDD} for the Caffarelli-Kohn-Nirenberg inequalities.

To be more precise, we are interested in showing that the  following trace inequality (see, e.g. \cite{BD}) can be improved for radially symmetric functions (in the first $n$ variables):

\begin{theorem*}
Let $u\in C_0^\infty(\mathbb{R}^n \times \mathbb{R}_+)$. Then, the following inequality holds
\begin{equation*}
\||x|^{-\beta} f(x,0) \|_{L^q(\mathbb{R}^n)} \le C \| |(y,z)|^{\alpha} \nabla f(y,z)\|_{L^2(\mathbb{R}^n\times \mathbb{R}^+)}
\end{equation*}provided that:
\begin{equation} 
0\le \alpha+\beta \le \frac12,
\end{equation}
\begin{equation}
\alpha> - \frac{n+1}{2}+1,
\end{equation}
and
\begin{equation}
\frac{n}{q}-\frac{n+1}{2}= \alpha+\beta -1.
\end{equation}
\end{theorem*}

Indeed, we will show that the  following refinement is possible in the case of radially symmetric functions:

\begin{theorem}
\label{teo-trazas}
Let $u\in C_0^\infty(\mathbb{R}^n \times \mathbb{R}_+)$ be a radially symmetric function in the first $n$ variables. Then, the following inequality holds
\begin{equation}
\label{desig-trazas}
\|f(x,0) |x|^{-\beta}\|_{L^q(\mathbb{R}^n)} \le C \| |(y,z)|^\alpha \nabla f(y,z)\|_{L^p(\mathbb{R}^n\times \mathbb{R}^+)}
\end{equation}
provided that:
\begin{equation} 
\label{alfamasbeta}
-\frac{n}{q'}\le \alpha +  \beta \le \frac{1}{p'},
\end{equation}
\begin{equation}
 \alpha > -\frac{n+1}{p}+1,
\end{equation}
and
\begin{equation}
\label{reescale-trazas}
\frac{n}{q}-\frac{n+1}{p}=\alpha + \beta -1.
\end{equation}
\end{theorem}

\begin{remark}
Using condition \eqref{reescale-trazas}, condition \eqref{alfamasbeta} can be seen to be equivalent to $1\le p \le q < \infty.$
\end{remark}

As a preliminary result for the proof of Theorem \ref{teo-trazas}, we will first prove the following theorem, of independent interest. 

\begin{theorem}
\label{operador-trazas}
Let $x\in \mathbb{R}^n$ and 
\begin{equation}
\label{op-trazas}
Tf(x):= \int_{\mathbb{R}^n \times \mathbb{R}_+} \frac{f(y,z)}{[(x-y)^2+z^2]^{\frac{n}{2}}} \, dy \, dz.
\end{equation}
Assume $f \in C_0^\infty(\mathbb{R}^n \times \mathbb{R}_+)$ is such that $f(y,z)= f_0(|y|,z)$. Then, the inequality 
\begin{equation}
\|Tf(x) |x|^{-\beta}\|_{L^q(\mathbb{R}^n)} \le C \| |(y,z)|^\alpha f(y,z) \|_{L^p(\mathbb{R}^n \times \mathbb{R}^+)} 
\end{equation}
holds provided that
\begin{equation}
\label{pq-traza}
1\le p\le q < \infty
\end{equation}
\begin{equation}
\label{reesc-op-trazas}
\frac{n}{q}-\frac{n+1}{p}=\alpha + \beta -1
\end{equation}
and
\begin{equation}
\label{beta-op-trazas}
-\frac{n}{q'}<\beta < \frac{n}{q}.
\end{equation}
\end{theorem}

The rest of this paper is organized as follows. In Section 2 we recall some necessary preliminaries. In Section 3 we prove Theorem \ref{ckn-radial}. In Section 4 we explain the relation between the operator $Tf$ defined by \eqref{op-trazas} and the weighted trace inequalities we are interested in, and find a convenient expression for this operator when acting on radially symmetric functions (in the first $n$ variables). In Section 5 we prove Theorem \ref{operador-trazas} and, finally, in Section 6, we use Theorem  \ref{operador-trazas}  to prove Theorem \ref{teo-trazas}.

\section{Notation and Preliminaries}

As it is usual, $C$ will denote a positive constant, independent of relevant parameters, that may change even within a single string of estimates. 

To prove Theorem \ref{ckn-radial} we will make use of a theorem proved in \cite{DDD}, that we recall here for the sake of completeness.

\begin{theorem}[\cite{DDD}, Theorem 1.2]
\label{teo-ddd}
For $n\ge 1$ define
\begin{equation}
\label{int-fr}(T_\gamma v)(x):= \int_{\R^n} \frac{v(y)}{|x-y|^\gamma} \; dy, \quad
0 <\gamma < n.
\end{equation}
Let
\begin{equation}
 1\le p\le q<\infty, 
\end{equation}
\begin{equation}
\alpha<\frac{n}{p'}
\end{equation}
\begin{equation}
 \beta<\frac{n}{q}
 \end{equation}
 \begin{equation}
 \label{ab}
\alpha + \beta \ge (n-1)(\frac{1}{q}-\frac{1}{p})
\end{equation}
and
\begin{equation}
\frac{1}{q}=\frac{1}{p}+\frac{\gamma+\alpha+\beta}{n}-1
\end{equation}
with strict inequality in \eqref{ab} if $p=1$. 
Then, the inequality
$$
\||x|^{-\beta}T_\gamma v\|_{L^q(\mathbb{R}^n)} \le C \||x|^\alpha
v\|_{L^p(\mathbb{R}^n)}
$$
holds for all radially symmetric $v\in L^p(\mathbb{R}^n, |x|^{p\alpha} dx)$,
where $C$ is independent of $v$.
\end{theorem}

For the proof of Theorem \ref{operador-trazas} we will use the main idea in the proof of Theorem \ref{teo-ddd}, that is, to write the operator (in this case, the operator given by \eqref{op-trazas} instead of the Riesz potential \eqref{int-fr}) as a convolution with respect to the Haar measure in $\mathbb{R}_+$. To make this precise, recall that if $G$ is a locally compact group, then $G$ posseses a Haar measure, that is, a
positive Borel measure $\mu$ which is left invariant (i.e., $\mu(At)=\mu(A)$
whenever $t\in G$ and $A\subseteq G$ is measurable)   and nonzero on nonempty
open sets. In particular, if  $G=\mathbb{R}-\{0\}$, then $\mu=\frac{dx}{|x|}$, and if $G=
\mathbb{R}_+$, then $\mu=\frac{dx}{x}$.

The convolution of two functions $f,g\in L^1(G)$  is defined as:
$$
(f*g)(x) = \int_G f(y) g(y^{-1}x) \, d\mu(y)
$$
where $y^{-1}$ denotes the inverse of $y$ in the group $G$.

The following version of Young's inequality holds in this setting:

\begin{theorem}\cite[Theorem 1.2.12]{G}
\label{young}
Let $G$ be a locally compact group with left Haar measure $\mu$ that satisfies
$\mu(A)=\mu(A^{-1})$ for all measurable $A\subseteq G$. Let
$1 \le p,q,s \le \infty$ satisfy
$$
\frac{1}{q}+1 = \frac{1}{p}+\frac{1}{s}.
$$
Then for all $f\in L^p(G, \mu)$ and $g\in L^s (G, \mu)$ we have
\be
\|f * g\|_{L^q(G, \mu)}\le \|g\|_{L^s (G, \mu)} \|f\|_{L^p(G,\mu)}.
\ee
\end{theorem}

\section{Proof of theorem \ref{ckn-radial}}

Clearly, when $a=0$ the theorem is completely trivial. Therefore, we will split the proof into two cases, namely, when $a=1$ and when $0<a<1$. 

\subsection{Case $a=1$} 

Notice  that in this case, $\sigma=\gamma$ by \eqref{sigma}.

Observing that for $u\in C_0^\infty(\R^n)$
$$
|u(x)| \le C \int_{\R^n} \frac{|\nabla u|(y)}{|x-y|^{n-1}} \, dy  := T_{n-1}(|\nabla u|)(x)
$$
we see that
 
 $$
 \| |x|^\gamma u\|_{L^r} \le C \| |x|^\gamma T_{n-1}(|\nabla u|)\|_{L^r}
 $$
but, since we are assuming that $u$ is a radial function, then so is $|\nabla u|$ and we can use  Theorem \ref{teo-ddd}
 to deduce that
$$
 \| |x|^\gamma T_{n-1}(|\nabla u|)\|_{L^r} \le C \| |x|^\alpha \nabla u\|_{L^p}
 $$
provided that
\begin{equation}
\label{pr}
1 \le p \le r <\infty 
\end{equation}
\begin{equation}
\label{reescale-pr}
\frac{1}{r} + \frac{\gamma}{n}= \frac{1}{p}+\frac{\alpha - 1}{n}
\end{equation}
\begin{equation}
\label{alfa}
\alpha < \frac{n}{p'}
\end{equation}
\begin{equation}
\label{gama-r}
-\gamma < \frac{n}{r}
\end{equation}
and
\begin{equation}
\label{a-g}
 (n-1)\left(\frac{1}{r}-\frac{1}{p} \right) \le \alpha -\gamma, 
\end{equation}
with strict inequality in \eqref{a-g} if $p=1$.

Clearly, the scaling condition \eqref{reescale-pr} equals condition \eqref{reescale} when $a=1$; and using \eqref{reescale-pr}, condition \eqref{pr} can be seen to be equivalent to $\gamma-\alpha \le 1$, which holds because of hypothesis \eqref{alfa-sigma} (recall that in this case $\gamma=\sigma$).  Condition \eqref{gama-r} equals condition \eqref{cond-sigma} (in this case it is  also included in \eqref{cond2}); and \eqref{a-g} follows from \eqref{alfa-sigma} since $a=1$.
 
We claim that condition \eqref{alfa} can be removed if we only wish to consider the inequality
\begin{equation}
\label{caso1}
\| |x|^\gamma u\|_{L^r} \le C \| |x|^\alpha \nabla u\|_{L^p}
\end{equation}
(this is not the case if the operator $T_{n-1}$ is not acting on $|\nabla u|$). Indeed, we will prove that if \eqref{caso1} holds for $\alpha$ and $\gamma$, then it also holds for $\alpha + 1$ and $\gamma + 1$, provided that  $\alpha p\neq -1$. To this end, we apply the inequality to $|x| u$ (strictly speaking, this function is not $C_0^\infty$, but it suffices to take a regularized distance function to the origin, see e.g. \cite{St}, and apply the same argument). 

Then,
$$
\||x|^{\gamma+1}u\|_r \le C \||x|^\alpha \nabla(|x|u)\|_p \sim C (\||x|^{\alpha+1} \nabla u\|_p + \||x|^\alpha u\|_p)
$$
and, therefore, it suffices to see that $\||x|^\alpha u\|_p \le C \||x|^{\alpha+1} \nabla u\|_p$. To this end write
\begin{align*}
\||x|^\alpha u\|_p^p &= \int |x|^{p\alpha} |u|^p \, dx 
\\ &\le C \int |\nabla |x|^{p\alpha +1}| |u|^p \, dx
\\ & \le C \int |x|^{p\alpha+1} |\nabla |u|^p| \, dx 
\\ &\le C \int |x|^{p\alpha+1} |u|^{p-1} |\nabla u| \, dx
\\ &\le C \left( \int  |x|^{p \alpha} |u|^p \, dx \right)^{\frac{1}{p'}} \left( \int |x|^{p(\alpha +1)} |\nabla u|^p  \, dx  \right)^{\frac{1}{p}}
\end{align*}

Thus, we have proved that
$$
\| |x|^\alpha u\|_p^p \le C \| |x|^\alpha u \|_p^{\frac{p}{p'}} \|  |x|^{\alpha+1} \nabla u\|_p
$$
whence it follows immediately that
$$
\||x|^\alpha u\|_p \le C \| |x|^{\alpha+1} \nabla u\|_p.
$$

Iterating the same argument, we can see that if \eqref{caso1} holds for  $\gamma$ and $\alpha$, then it also holds for  $\gamma+k$ and $\alpha+k$ with $k\in \mathbb{N}_0$ provided that $(\alpha -k)p \neq -1$. Therefore, to see that  we can remove condition \eqref{alfa}, it suffices to observe that any $\alpha \ge \frac{n}{p'}$ can be written as $(\alpha -k) + k$, with $-\frac{n}{p} < \alpha -k < \frac{n}{p'}$, and $(\alpha -k)p \neq -1$.  Indeed, since $\frac{n}{p'}-  (-\frac{n}{p})= n$, such a  $k$ exists except when $n=1$ and $\alpha=\frac{1}{p'}$. But this is impossible, since in that case, by \eqref{reescale-pr}
 we should have $\frac{1}{r}+\gamma = \frac{1}{p}+\frac{1}{p'}-1$, that is, $ \frac{1}{r}+\gamma=0$, which contradicts \eqref{cond2}.

\subsection{Case $0<a<1$} 

Write
\begin{align}
\nonumber \left( \int |x|^{\gamma r} |u|^{ra+(1-a)r} \, dx  \right)^{\frac{1}{r}} & = \left( \int |x|^{r \beta (1-a)} |u|^{(1-a)r} |x|^{r\gamma (1-\frac{\beta (1-a)}{ \gamma})}  |u|^{ar} \, dx \right)^{\frac{1}{r}}
\\ \label{holder} &\le  \| |x|^\beta u \|_{L^q}^{1-a} \| |x|^{\frac{\gamma}{a} (1-\frac{\beta(1-a)}{\gamma})} u\|_{L^{\frac{arq}{q-r(1-a)}}}^a
\\ \label{uso-sigma} & =  \| |x|^\beta u \|_{L^q}^{1-a} \||x|^\sigma u\|_{L^\frac{arq}{q-r(1-a)}}^a
\end{align}
where in \eqref{holder}  we have used H\"older's inequality with exponent $\frac{q}{r(1-a)}$ (which is greater than 1 by \eqref{cond-qr}) and in \eqref{uso-sigma} we have used the definition of $\sigma$, given in \eqref{sigma}.

Applying now the result obtained in the case $a=1$, we deduce that 
\begin{equation*}
\| |x|^\gamma u\|_{L^r} \le C \| |x|^\beta u\|_{L^q}^{1-a} \| |x|^\alpha \nabla u\|_{L^p}^a
\end{equation*}
provided that
\begin{equation}
\label{gen4}
1\le p\le \frac{arq}{q-r(1-a)} < \infty
\end{equation}
\begin{equation}
\label{gen3}
\frac{q-r(1-a)}{arq} + \frac{\sigma}{n} = \frac{1}{p} + \frac{\alpha -1}{n}
\end{equation}
\begin{equation}
\label{gen1}
 -\sigma < \frac{n(q-r+ar)}{arq} 
\end{equation}
and
\begin{equation}
\label{gen2}
\alpha- \sigma \ge (n-1)\left( \frac{q -r(1-a)}{arq} -\frac{1}{p}\right),
\end{equation}
where in \eqref{gen2} the inequality is strict if $p=1$.

Clearly,  condition \eqref{gen4} holds because of \eqref{cond-qr}, and condition \eqref{gen3} is easily seen to be equivalent to \eqref{reescale} using the definition of $\sigma$ given in \eqref{sigma}. Finally, condition \eqref{gen1} equals conditon \eqref{cond-sigma} while \eqref{gen2} is the same as \eqref{alfa-sigma}. This concludes the proof.

\section{The operator associated to trace inequalities}

Before we can proceed to the proof of the announced trace inequality, we first need to obtain an expression analogous to \eqref{representation} and, then, a convenient expression for the involved operator when acting on radial functions.

To this end, given  $u$ and a unitary vector $\xi$, consider $g(s)=u(s\xi, 0)$. Then, $g(0)= - \int_0^\infty g'(s) \, ds = - \int_0^\infty \nabla u (s \xi) \cdot \xi \, ds$. 

Consider now $\varphi \in C_0^\infty(S^n)$ supported in $\mathbb{R}^n \times \mathbb{R}_+$ and such that $\int_{Sn} \varphi(\xi) \, d\sigma (\xi)=1$. Then
$$
u(0,0)= - \int_0^\infty \int_{S_n} \nabla u (s\xi) \cdot \xi \,  \varphi(\xi) \, d\xi \, ds.
$$

For $(y,z)\in \mathbb{R}^{n+1}$ let $\phi(y,z)=\varphi((y,z) / \|(y,z)\|)$. Therefore,  $\phi(s\xi) = \varphi(\xi)$ for all $s\in \mathbb{R}^+, \xi \in S^n$, and the above identity becomes
\begin{align*}
u(0,0) &= - \int_0^\infty \int_{S_n} \nabla u (s\xi) \cdot s\xi \, \phi(s\xi) \frac{1}{s^{n+1}} s^n \, ds \, d\xi
\\ & = -\int_{\mathbb{R}^n \times \mathbb{R}_+} \nabla u (y,z) \cdot (y,z) \, \phi(y,z) \frac{1}{\|(y,z)\|^{n+1}} \, dy \, dz 
\end{align*}

More generally,
\begin{align*}
|u(x,0)| & \le \int_{\mathbb{R}^n \times \mathbb{R}_+} |\nabla u (y,z)| \frac{1}{\|(x-y,z)\|^n} \, dy \, dz
\\ & = \int_{\mathbb{R}^n \times \mathbb{R}_+} |\nabla u(y,z)| \frac{1}{[(x-y)^2+z^2]^{\frac{n}{2}}} \, dy \, dz
\end{align*}

Then, we have to study the behavior of the operator
$$
Tf(x)= \int_{\mathbb{R}^n \times \mathbb{R}_+} \frac{f(y,z)}{[(x-y)^2+z^2]^{\frac{n}{2}}} \, dy \, dz
$$
for $x\in \mathbb{R}^n$.

Since we are interested in the radial case, assume $f$ is a radially symmetric function in the first variable (by an abuse of notation we will still call it  $f$).

Using polar coordinates
\begin{align*}
&y =ry', \quad r=|y|,   \quad y'\in S^{n-1}
\\ & x=\rho x' , \quad \rho=|x|, \quad x'\in S^{n-1}
\end{align*}
if $n\ge 2$ we may write:
\begin{align*}
Tf(x) &= \int_0^\infty \left[ \int_0^\infty \int_{S^{n-1}} \frac{f(r,z) \, r^{n-1} }{(\rho^2 - 2 \rho r x' \cdot y' + r^2 + z^2)^{\frac{n}{2}} } \, dy' \, dr \right] \, dz
\\ &= \int_0^\infty  \int_0^\infty f(r,z) r^{n-1}  \int_{-1}^1 \frac{(1-t^2)^{\frac{n-3}{2}}}{(\rho^2 - 2 \rho r t + r^2 + z^2)^{\frac{n}{2}}} \, dt \, dr  \, dz
\end{align*}
where the second equality can be justified integrating in the sphere (see, e.g., Lemma 4.1 from \cite{DDD}).

Making the change of variables $ z = r \bar z$, $dz = r \, d\bar z$ we obtain
\begin{align}
\nonumber Tf(x) &= \int_0^\infty  \int_0^\infty f(r,r\bar z) r^n  \int_{-1}^1 \frac{(1-t^2)^{\frac{n-3}{2}}}{ r^n \left[1 - 2 \left(\frac{\rho}{r}\right) t + \left(\frac{\rho}{r}\right)^2 +  \bar z^2\right]^{\frac{n}{2}}} \, dt \, dr  \, d\bar z
\\ &= \int_0^\infty \int_0^\infty f(r, rz)  I\left(\frac{\rho}{r}, z\right) \, dr \, dz \label{expre}
\end{align}
where, for $a>0$, 
$$
I(a,z) :=  \int_{-1}^1 \frac{(1-t^2)^{\frac{n-3}{2}}}{(1 -2at + a^2 + z^2)^{\frac{n}{2}}} \, dt.
$$
Expression \eqref{expre} will allow us to write $Tf$ as convolution operator and to obtain Theorem \ref{operador-trazas}, that we proceed to prove next.

\section{Proof of Theorem  \ref{operador-trazas}}

If $n=1$ recall that we want to prove
$$
\|Tf(x) |x|^{-\beta}\|_{L^q(\mathbb{R})} \le C \| |(y,z)|^\alpha f(y,z) \|_{L^p(\mathbb{R} \times \mathbb{R}^+)} 
$$
Since in this case \eqref{expre} does not hold, we remark that
$$
\|Tf(x) |x|^{-\beta}\|_{L^q(\mathbb{R}, dx)} = \| |x|^{-\beta+\frac{1}{q}} Tf \|_{L^q(\mathbb{R}, \frac{dx}{|x|})}
$$
and write
\begin{align*}
|x|^{-\beta + \frac{1}{q}} Tf(x) &= \int_{-\infty}^\infty \int_0^\infty \frac{f(y,z) |x|^{-\beta+\frac{1}{q}}}{[(x-y)^2 + z^2]^{\frac12}} \, dz \,dy
\\ & = \int_{-\infty}^\infty \int_0^\infty \frac{f(y,|y| \bar z) |x|^{-\beta+\frac{1}{q}} |y|}{( | \frac{x}{y}-1|^2 + \bar z^2)^{\frac12}} \, d\bar z \, \frac{dy}{|y|}
\\ & = \int_0^\infty \int_{-\infty}^\infty  \frac{f(y,|y| \bar z) (\frac{|x|}{|y|})^{-\beta+\frac{1}{q}} |y|^{1-\beta+\frac{1}{q}}}{( | \frac{x}{y}-1|^2 + \bar z^2)^{\frac12}}\, \frac{dy}{|y|}  \, d\bar z 
\\ & = \int_0^\infty (f(y, |y|\bar z) |y|^{1-\beta+\frac{1}{q}}) * \left( \frac{|y|^{-\beta+\frac{1}{q}}}{ (|y-1|^2 + \bar z^2)^\frac12 }\right) \, d\bar z
\end{align*}
where the convolution is taken with respect to the first variable in the multiplicative group $\mathbb{R}-\{0\}$ with Haar measure $dx/|x|$.

Let  $g(y)=\frac{|y|^{-\beta+\frac{1}{q}}}{ (|y-1|^2 + \bar z^2)^\frac12}$. Then, by Young's inequality (Theorem \ref{young}), if 
\begin{equation}
\label{you}
\frac{1}{q}+1 = \frac{1}{p} + \frac{1}{s}
\end{equation}
\begin{align}
\nonumber \|Tf(x) & |x|^{-\beta}\|_{L^q(\mathbb{R}, dx)}  
\\ \nonumber & \le \int_0^\infty \| f(y,|y|\bar z) |y|^{1-\beta+\frac{1}{q}} \|_{L^p(\frac{dy}{|y|})} \| g  \|_{L^s(\frac{dy}{|y|})} \, d\bar z
\\ & = \int_0^\infty \left( \int_{-\infty}^\infty |f(y,|y| \bar z)|^p |y|^{(1-\beta+\frac{1}{q})p -1} (1+\bar z^2)^{\frac{\alpha p}{2}} \right)^{\frac{1}{p}} \left( \frac{\|g\|_{L^s(\frac{dy}{|y|})}}{(1+\bar z^2)^{\frac{\alpha}{2}}}\right)  \, d\bar z \label{ult}
\end{align}

Observing now that 
\begin{align*}
\| |(y,z)|^\alpha f(y,z) \|_{L^p(\mathbb{R}\times \mathbb{R}_+)} & = \int_0^\infty \int_{-\infty}^\infty (y^2+z^2)^{\frac{\alpha p}{2}} |f(y,z)|^p \, dy \, dz
\\ & = \int_0^\infty \int_{-\infty}^\infty (y^2 + y^2 \bar z^2)^{\frac{\alpha p}{2}} |f(y, |y|\bar z)|^p |y| \, dy \, d\bar z
\\ & = \int_0^\infty \int_{-\infty}^\infty (1+\bar z^2)^{\frac{\alpha p}{2}} |y|^{\alpha p +1} |f(y, |y|\bar z)|^p \, dy \, d\bar z
\end{align*}
and that $(1-\beta+\frac{1}{q}) p -1 = \alpha p +1$ (by \eqref{reesc-op-trazas}), we can apply H\"older's inequality to \eqref{ult} to obtain
$$
\|Tf(x) |x|^{-\beta}\|_{L^q(\mathbb{R}^n)} \le  \| |(y,z)|^\alpha f(y,z) \|_{L^p(\mathbb{R}^n \times \mathbb{R}^+)}  \left( \int_0^\infty \frac{\|g\|_{L^s(\frac{dy}{|y|})}^{p'}}{(1+z^2)^{\frac{\alpha p'}{2}}} \, d z \right)^{\frac{1}{p'}}
$$
Therefore, to conclude the proof of the one-dimensional case it suffices to see that
$$
\int_0^\infty \frac{\|g\|_{L^s(\frac{dy}{|y|})}^{p'}}{(1+z^2)^{\frac{\alpha p'}{2}}} \, dz  < + \infty
$$
provided that \eqref{reesc-op-trazas}, \eqref{beta-op-trazas}  and \eqref{you} hold. We omit the details since the computations are analogous to those that we will do in the higher dimensional case.

Now we proceed to the case $n \ge 2$. In this case, remark that, 
\begin{align*}
\|Tf(x) |x|^{-\beta}\|_{L^q(\mathbb{R}^n)} & = C \left( \int_0^\infty |Tf(\rho)|^q \rho^{-\beta q+n} \frac{d\rho}{\rho}\right)^{\frac{1}{q}} 
\\ &=  C \|  \rho^{-\beta+\frac{n}{q}} Tf\|_{L^q(\frac{d\rho}{\rho})}
\end{align*}

We claim that $\rho^{-\beta+\frac{n}{q}} Tf$ can be written as a convolution in the multiplicative group $(\mathbb{R}_+, \cdot)$. Indeed,
\begin{align*}
\rho^{-\beta+\frac{n}{q}} Tf & = \int_0^\infty \int_0^\infty f(r,rz)  \, I\left(\frac{\rho}{r}, z\right) \rho^{-\beta+\frac{n}{q}} \, dr dz
\\ &= \int_0^\infty \int_0^\infty f(r,rz) \, I\left(\frac{\rho}{r}, z \right) \left(\frac{\rho}{r} \right)^{-\beta+\frac{n}{q}} r^{-\beta+\frac{n}{q}+1} \, \frac{dr}{r} dz
\\ &= \int_0^\infty (f(r,rz) r^{-\beta + \frac{n}{q}+1}) * (I(r,z) r^{-\beta+\frac{n}{q}}) \, dz
\end{align*}
where $*$ denotes the convolution with respect to the Haar measure $dr/r$ in the first variable.

Therefore, using Young's inequality, for
\begin{equation}
\label{young-coef}
\frac{1}{q} +1 =\frac{1}{p}+\frac{1}{s},
\end{equation}
 we obtain 
\begin{align}
 \nonumber \| Tf(&\rho) \rho^{-\beta + \frac{n}{q}}\|_{L^q(\frac{d\rho}{\rho})}
\\ \nonumber & \le \int_0^\infty \|(f(r,rz) r^{-\beta + \frac{n}{q}+1}) * (I(r,z) r^{-\beta+\frac{n}{q}})\|_{L^q(\frac{dr}{r})} \, dz
\\ \nonumber & \le \int_0^\infty \|f(r,rz) r^{-\beta+\frac{n}{q}+1}\|_{L^p(\frac{dr}{r})} \|I(r,z) r^{-\beta+\frac{n}{q}}\|_{L^{s}(\frac{dr}{r})} \, dz
\\ \nonumber &= \int_0^\infty \left( \int_0^\infty |f(r,rz)|^p r^{(-\beta+\frac{n}{q}+1)p} \frac{dr}{r} \right)^{\frac{1}{p}}  \|I(r,z) r^{-\beta+\frac{n}{q}}\|_{L^{s}(\frac{dr}{r})} \, dz
\\ \nonumber & = \int_0^\infty \left( \int_0^\infty |f(r,rz)|^p r^{(-\beta+\frac{n}{q}+1)p} (1+z^2)^{\frac{\alpha p}{2}} \frac{dr}{r} \right)^{\frac{1}{p}}  \frac{\|I(r,z) r^{-\beta+\frac{n}{q}}\|_{L^{s}(\frac{dr}{r})}}{(1+z^2)^{\frac{\alpha}{2}}}   \, dz 
\end{align}

Now, since \begin{align*}
\| |(y,z)|^\alpha f(y,z) \|_{L^p(\mathbb{R}^n \times \mathbb{R}^+)} &= \left( \int_0^\infty \int_0^\infty (r^2+z^2)^{\frac{\alpha p}{2}} |f(r,z)|^p r^{n-1} \, dr dz \right)^{\frac{1}{p}}
\\ &= \left( \int_0^\infty \int_0^\infty (r^2+r^2\bar z^2)^{\frac{\alpha p}{2}} |f(r,r\bar z)|^p r^n \, d\bar z dr \right)^{\frac{1}{p}}
\\&= \left( \int_0^\infty \int_0^\infty r^{\alpha p} (1+\bar z^2)^{\frac{\alpha p}{2}} |f(r,r \bar z)|^p r^n \, d \bar z dr \right)^{\frac{1}{p}},
\end{align*}
observing that $n+\alpha p = p (-\beta + \frac{n}{q} + 1) -1$ and applying H\"older's inequality to the above  expression, we obtain
\begin{align*}
\| Tf &(\rho) \rho^{-\beta + \frac{n}{q}}\|_{L^q(\frac{d\rho}{\rho})} 
\\ &\le  \| |(y,z)|^\alpha f(y,z) \|_{L^p(\mathbb{R}^n \times \mathbb{R}^+)} \left( \int_0^\infty   \frac{\|I(r,z) r^{-\beta+\frac{n}{q}}\|^{p'}_{L^{s}(\frac{dr}{r})}}{(1+z^2)^{\frac{\alpha p'}{2}}}   \, dz \right)^{\frac{1}{p'}}
\end{align*}

Therefore, to conclude the proof of the theorem it suffices to see that
\begin{equation}
\label{norma}
\int_0^\infty   \|I(r,z) r^{-\beta+\frac{n}{q}}\|^{p'}_{L^{s}(\frac{dr}{r})} (1+z^2)^{-\frac{\alpha p'}{2}}  \, dz < +\infty.
\end{equation}

Observe first that the denominator of
$$
I(r,z)=  \int_{-1}^1 \frac{(1-t^2)^{\frac{n-3}{2}}}{(1 -2rt + r^2 + z^2)^{\frac{n}{2}}} \, dt
$$
can be rewritten as  $[(a-t)^2+(1-t^2) + z^2]^{\frac{n}{2}}$ and, therefore, it vanishes for $r=t=1$ and $z=0$ only.

To bound $ \|I(r,z) r^{-\beta+\frac{n}{q}}\|_{L^{s}(\frac{dr}{r})}$, consider $\varphi \in C^\infty(\mathbb{R})$ such that $supp( \varphi) \subseteq [\frac12, \frac32]$, $0\le \varphi \le 1$ and  $\varphi \equiv 1$ in $(\frac34,\frac54)$. We can then split $I(r,z) r^{-\beta+\frac{n}{q}}= I(r,z) r^{-\beta+\frac{n}{q}} \varphi(r) + I(r,z) r^{-\beta+\frac{n}{q}} (1-\varphi(r)) = g_1(r) + g_2(r)$ and bound both terms separately. To this end, we will study first the behavior of $g_1$ and $g_2$ and then estimate \eqref{norma}.

Consider first  $g_2$. For $r\to 0$, we have
 $$I(0,z) = (1+z^2)^{-\frac{n}{2}} \int_{-1}^1 (1-t^2)^{\frac{n-3}{2}} \, dt \sim (1+z^2)^{-\frac{n}{2}}.$$
 Therefore,  $\|g_2\|_{L^s(\frac{dr}{r})}$ behaves like $(1+z^2)^{-\frac{n}{2}} $ when $r\to 0$, provided that $\beta < \frac{n}{q}$.

When $r\to \infty$,
$$
I(r,z)\sim \frac{1}{(r^2 + z^2)^{\frac{n}{2}}}.
$$
In this case, if $z$ is bounded,  say $z\le 2$, $\|g_2\|_{L^s(\frac{dr}{r})}$ is also bounded provided that $\beta>-\frac{n}{q'}$. On the other hand, when $z\to \infty$, we need to estimate
\begin{align*}
\left( \int_2^\infty \frac{r^{s(-\beta+\frac{n}{q})}}{(r^2+z^2)^{\frac{ns}{2}}} \, \frac{dr}{r} \right)^{\frac{1}{s}} &=  \left( z^{s(-\beta+\frac{n}{q}-n)} \int_{\frac{2}{z}}^\infty \frac{r^{s(-\beta+\frac{n}{q})}}{(r^2+1)^{\frac{ns}{2}}} \, \frac{dr}{r}\right)^{\frac{1}{s}} 
\\ &\sim z^{-\beta-\frac{n}{q'}}
\end{align*}
assuming again that $\beta>-\frac{n}{q'}$.

We can proceed now to $\|g_1\|_{L^s(\frac{dr}{r})}$.  We consider first the case  $k=\frac{n-3}{2} \in \N_0$, that is $n\ge 3$ and odd. 

If $z$ is sufficiently large, then $I(r,z) \sim z^{-n}$ and, therefore, $\|g_1\|_{L^s(\frac{dr}{r})} \sim z^{-n}$.

If, on the contrary, $z\to 0$, we may write
$$
I(r,z) \sim \int_{-1}^1 (1-t^2)^k \frac{d^k}{dt^k} \,  \left\{(1-2rt+r^2 +z^2)^{-\frac{n}{2}+k} \right\} \, dt
$$
and integrating  by parts $k$-times  (the boundary terms vanish), we obtain 
$$I(r,z) \le C_k [(1-r)^2 + z^2]^{-\frac{n}{2}+k+1}.$$
Since we are assuming that $-\frac{n}{2}+k+1=-\frac{1}{2}$, we conclude that
\begin{align*}
\|g_1\|_{L^s(\frac{dr}{r})}            & \sim \left(  \int_{\frac12}^{\frac32} \frac{dr}{[(1-r)^2 + z^2]^{\frac{s}{2}}}\right)^{\frac{1}{s}} 
\\ &\sim \left(  \int_{\frac12}^{\frac32} \frac{dr}{(|1-r|+ z)^s}\right)^{\frac{1}{s}} 
\\ & \sim \frac{1}{z^{1-\frac{1}{s}}}
\end{align*}

We can consider now  $k=m+\frac12, m \in \mathbb{N}_0$. In this case
\begin{align*}
|\frac{d}{dz} & I(r,z) |
\\ &\le C z \int_{-1}^1 \frac{ (1-t^2)^k}{(1-2rt+r^2+z^2)^{\frac{n}{2}+1}} \, dt
\\ & \leq C z \left( \int_{-1}^1 \frac{(1-t^2)^m}{(1-2rt+r^2 +z^2)^{\frac{n+2}{2}}} \, dt\right)^{\frac12}  \left( \int_{-1}^1 \frac{(1-t^2)^{m+1}}{(1-2rt+r^2 +z^2)^{\frac{n+2}{2}}} \, dt\right)^{\frac12}
\end{align*}
and, since now $\frac{n+2}{2} \in \N$, we deduce from the previous case that
\begin{align*}
\left|\frac{d}{dz}  I(r,z) \right| &\le C z [(1-r^2) + z^2]^{\frac{-(n+2)+2m+3}{2}} 
\\ &=  C z [(1-r)^2 + z^2]^{-\frac32} 
\\ &\le C z  [|1-r| + z]^{-3} 
\end{align*}

Therefore, 
$$
I(r,z) = \int_0^z  \frac{d}{dt}I (r,t) \, dt         \le C z [|1-r| + s]^{-2} |_0^z \le C z [|1-r| + z]^{-2} 
$$
which implies
$$
\|g_1\|_{L^s(\frac{dr}{r})} \sim \frac{1}{z^{1-\frac{1}{s}}}.
$$

It remains to check the case $k= -\frac12$ (i.e., $n=2$). To this end, we write
\begin{align*}
I(r,z) &= \underbrace{ \int_{-1}^0 \frac{(1-t^2)^{-\frac12}}{(1-2at+a^2+z^2)} \, dt}_{(i)} + \underbrace{ \int_0^1 \frac{(1-t^2)^{-\frac12}}{(1-2at+a^2+z^2)} \, dt }_{(ii)}
\end{align*}
Clearly, 
$$
(i) \le \int_{-1}^0 \frac{ dt}{(1+t)^{\frac12}} = 2
$$
while
\begin{align*}
(ii) &\le \int_0^1 \frac{(1-t)^{-\frac12}}{(1-2at+a^2+z^2)} \, dt 
\\ &= -2 \int_0^1 \frac{ \frac{d}{dt}[(1-t)^{\frac12}]}{1-2at+a^2+z^2} \, dt
\\ &\le 4a \int_0^1 \frac{(1-t^2)^{\frac12}}{(1-2at+a^2+z^2)^2} \, dt
\end{align*}
and the last integral can be bounded as before (notice that it corresponds to the case $n=4$).

We are now able to see that \eqref{norma} holds. Indeed, by our previous calculations, we need to bound
\begin{align*}
& \int_0^1 \left( \frac{1}{z^{1-\frac{1}{s}} (1+z^2)^{\frac{\alpha}{2}}} + \frac{1}{(1+z^2)^{\frac{n+\alpha}{2}}}  \right)^{p'}  dz 
\\ + &   \int_1^\infty  \left( \frac{1}{z^n (1+z^2)^{\frac{\alpha}{2}}} + \frac{1}{z^{\beta+\frac{n}{q'}}(1+z^2)^{\frac{\alpha}{2}}}  \right)^{p'}  dz
\end{align*}

When $z\to 0$, the integrability condition is $p'(1-\frac{1}{s})<1$, which holds because of \eqref{pq-traza} and \eqref{young-coef}. When $z\to \infty$, since we are assuming that  $\beta<\frac{n}{q}$, there holds that $n>\beta+\frac{n}{q'}$, whence the integralibity condition is  $p'(\beta+\frac{n}{q'}+\alpha)>1$, that is, $\alpha + \beta>\frac{1}{p'}-\frac{n}{q'}$. But, by \eqref{reescale-trazas} this condition is equivalent to $\frac{n}{p'}>0$, which trivially holds. This concludes the proof of the theorem.

\section{Proof of Theorem \ref{teo-trazas}}

As in the case of the Caffarelli-Kohn-Nirenberg interpolation inequality, if we simply apply Theorem \ref{operador-trazas} 
 to $|\nabla f|$ we obtain \eqref{desig-trazas}
provided that
\begin{equation}
1\le p\le q < \infty
\end{equation} 
\begin{equation}
\label{reesc-traza}
 \frac{n}{q}-\frac{n+1}{p}=\alpha + \beta -1
 \end{equation}
 and
 \begin{equation}
 \label{beta}
-\frac{n}{q'}<\beta < \frac{n}{q}.
 \end{equation}
 Notice that this last condition is equivalent to $-\frac{n+1}{p}+1 <\alpha < \frac{n+1}{p'}$ because of \eqref{reesc-traza}.

To prove Theorem  \ref{teo-trazas} we need to see that  condition $\alpha<\frac{n+1}{p'}$ is unnecessary for inequality \eqref{desig-trazas} to hold. Indeed, with a similar argument as that used for Theorem \ref{ckn-radial}, we will prove that if the inequality holds for  $\alpha$ and $\beta$ then it also holds for  $\alpha+1$ and $\beta-1$ provided that $\alpha p \neq -1$.

To see this, consider $f(x)|x|$ (strictly speaking, we would need to replace $|x|$ by a regularized distance, to guarantee that the product is in $C_0^\infty$). Then, 
\begin{align*}
\|f(x,0) |x|^{-\beta+1}\|_{L^q(\mathbb{R}^n)} &\le C \| |(y,z)|^\alpha \nabla(|(y,z)| f(y,z))\|_{L^p(\mathbb{R}^n\times \mathbb{R}^+)}
\\ &\le C \||(y,z)|^{\alpha+1} \nabla f(y,z)\|_{L^p(\mathbb{R}^n\times \mathbb{R}^+)} 
\\  & \quad + \||(y,z)|^\alpha f(y,z)\|_{L^p(\mathbb{R}^n\times \mathbb{R}^+)}
\end{align*}

Therefore, it suffices to see that 
$$
 \||(y,z)|^\alpha f(y,z)\|_{L^p(\mathbb{R}^n\times \mathbb{R}^+)} \le C \||(y,z)|^{\alpha+1} \nabla f(y,z)\|_{L^p(\mathbb{R}^n\times \mathbb{R}^+)}.
 $$ 

To this end, consider
\begin{align*}
\|| (y,z)|^\alpha & f(y,z)\|_{L^p(\mathbb{R}^n\times \mathbb{R}^+)}^p 
\\ & = \int_{\mathbb{R}_+} \int_{\mathbb{R}^n} |(y,z)|^{p\alpha} |f(y,z)|^p \, dy \, dz 
\\ &\le C \int_{\mathbb{R}_+} \int_{\mathbb{R}^n}  |\nabla |(y,z)|^{p\alpha +1}| |f(y,z)|^p \, dy \, dz
\\ &\le C \int_{\mathbb{R}_+} \int_{\mathbb{R}^n}  |(y,z)|^{p\alpha+1} |\nabla |f(y,z)|^p| \, dy \, dz
\\ & \le C \int_{\mathbb{R}_+} \int_{\mathbb{R}^n}  |(y,z)|^{p\alpha+1} |f(y,z)|^{p-1} |\nabla f(y,z)| \, dy \, dz
\\ &= C \int_{\mathbb{R}_+} \int_{\mathbb{R}^n}  |(y,z)|^{\alpha(p-1)}  |f(y,z)|^{p-1} |(y,z)|^{\alpha +1} |\nabla f(y,z)|  \, dy \, dz
\end{align*}

Applying H\"older's inequality we see that
$$
\| |(y,z)|^\alpha f(y,z)\|_p^p \le C \| |(y,z)|^\alpha f(y,z) \|_p^{\frac{p}{p'}} \|  |(y,z)|^{\alpha+1} \nabla f(y,z)\|_p
$$
and it follows immediately that
$$
\||(y,z)|^\alpha f(y,z)\|_p \le C \| |(y,z)|^{\alpha+1} \nabla f(y,z)\|_p
$$
as we wanted to see.

Iterating the same argument we see that if inequality \eqref{desig-trazas} holds for  $\alpha$ and $\beta$, then it holds for  $\alpha+k$ and $\beta-k$ with $k\in \mathbb{N}_0$. Therefore, to see that condition $\alpha<\frac{n+1}{p'}$ is uneccessary, it suffices to see that any $\alpha \ge \frac{n+1}{p'}$ can be written as $(\alpha -k) +k $, with $-\frac{n+1}{p}+1 <\alpha -k < \frac{n+1}{p'}$ and $(\alpha -k)p \neq -1$. 

But,  $\frac{n+1}{p'}- \left( -\frac{n+1}{p} +1\right)= n$, and therefore  $k$ can be chosen as above, except when  $n=1$ and $\alpha=\frac{n+1}{p'} = \frac{2}{p'}$ (that is, $\beta=-\frac{1}{q'}$) that cannot happen because for $n=1$, $\alpha>\frac{2}{p'}$ (because of \eqref{beta} and \eqref{reesc-traza}). This completes the proof of the theorem.

\end{document}